\documentclass[12pt]{article}

\usepackage{geometry}
\usepackage{graphicx}
\usepackage{enumerate}
\usepackage{subcaption}
\usepackage{adjustbox}

\usepackage[hidelinks]{hyperref}
\usepackage{breakurl}
\usepackage{url}
\usepackage{xcolor}

\usepackage{amssymb,amsthm,amsmath,amsfonts,mathtools}
\newtheorem{theorem}{Theorem}

\newcommand{\secref}[1]{\S\ref{#1}}

\newcommand{\ie}{\emph{i.e.}}
\newcommand{\eg}{\emph{e.g.}}

\newcommand{\RR}{\mbox{\bf R}}
\newcommand{\ZZ}{\mbox{\bf Z}}
\newcommand{\reals}{{\mbox{\bf R}}}
\newcommand{\Set}[2]{\{\,#1\,\mid\,#2\,\}}

\DeclarePairedDelimiter\ceil{\lceil}{\rceil}
\DeclarePairedDelimiter\floor{\lfloor}{\rfloor}
\DeclarePairedDelimiter{\norm}{\|}{\|}

\newcommand{\LEAF}{\mathrm{LEAF}}

\newcommand{\AFF}{\mathrm{AFF}}
\newcommand{\CVX}{\mathrm{CVX}}
\newcommand{\CCV}{\mathrm{CCV}}
\newcommand{\QCVX}{\mathrm{QCVX}}
\newcommand{\QCCV}{\mathrm{QCCV}}

\newcommand{\lmbdmax}{\lambda_{\mathrm{max}}}

\usepackage{listings}
\usepackage[
    natbib=true,
    citestyle=alphabetic,
    style=alphabetic,
    maxnames=3,
    minnames=3,
    maxbibnames=99,
    backend=bibtex,
    urldate=iso8601,
    firstinits=true,
    isbn=false,
    doi=false,
    date=year,
    sorting=nty
]{biblatex}
\DeclareFieldFormat[article,inbook,incollection,inproceedings,patent,thesis,
  unpublished]{citetitle}{#1}
\DeclareFieldFormat[article,inbook,incollection,inproceedings,patent,thesis,
  unpublished]{title}{#1}
\renewbibmacro{in:}{%
  \ifentrytype{article}{}{\printtext{\bibstring{in}\space}}}
\renewcommand*{\labelalphaothers}{\textsuperscript{+}}

\title{Disciplined Quasiconvex Programming}
\author{
  Akshay Agrawal \\ \texttt{
  \small akshayka@cs.stanford.edu} \and
  Stephen Boyd \\ \texttt{\small boyd@stanford.edu}
}

\begin{document}
\maketitle

\begin{abstract}
We present a composition rule involving quasiconvex functions that
generalizes the classical composition rule for convex functions. This rule
complements well-known rules for the curvature of quasiconvex functions under
increasing functions and pointwise maximums. We refer to the class of
optimization problems generated by these rules, along with a base set of
quasiconvex and quasiconcave functions, as \emph{disciplined quasiconvex
programs}. Disciplined quasiconvex programming generalizes
disciplined convex programming, the class of optimization problems targeted by
most modern domain-specific languages for convex optimization.  We describe an
implementation of disciplined quasiconvex programming that makes it possible to
specify and solve quasiconvex programs in CVXPY 1.0.
\end{abstract}

\section{Introduction}
A real-valued function $f$ is \emph{quasiconvex} if its domain
$C$ is convex, and for any $\alpha\in \reals$, its
$\alpha$-sublevel sets $\Set{x \in C}{f(x) \leq \alpha}$
are convex \cite[\S 3.4]{boyd2004}.
A function $f$ is quasiconcave if $-f$ is quasiconvex, and it
is quasilinear if it is both quasiconvex and quasiconcave. A
\emph{quasiconvex program} (QCP) is a mathematical optimization problem in
which the objective is to minimize a quasiconvex function over a convex set.
Because every convex function is also quasiconvex, QCPs generalize convex
programs. Though QCPs are in general nonconvex, many can nonetheless be solved
efficiently by a bisection method that involves
solving a sequence of convex programs \cite[\S 4.2.5]{boyd2004}, or by
subgradient methods \cite{kiwiel2001, konnov2003}.

The study of quasiconvex functions is several decades old \cite{fenchel53, nikaido54,
luenberger68}. Quasiconvexity has been of particular interest in
economics, where it arose in the study of competitive equilibria and the
modeling of utility functions \cite{arrow1954, guerraggio2004}. More recently,
quasiconvex programming has been applied to control \cite{gu1994, bullo2006,
seiler2010}, model order reduction \cite{sou2008}, computer vision
\cite{ke2006, ke2007}, computational geometry \cite{eppstein2005}, and machine
learning \cite{hazan2015}. While QCPs have many applications, it remains
difficult for non-experts to specify and solve them in practice. The point of
this paper is to close that gap.

Domain-specific languages (DSLs) have made convex optimization widely
accessible. DSLs let users specify their programs in natural mathematical
notation, abstracting away the process of canonicalizing problems to standard
forms for numerical solvers. The syntax of most DSLs for convex
optimization, including CVX \citep{cvx}, CVXPY \citep{cvxpy, cvxpyrewriting},
Convex.jl \citep{convexjl}, and CVXR \citep{fu2017cvxr}, is determined by a
grammar known as \emph{disciplined convex programming} (DCP)
\citep{grant2006disciplined}. DCP includes a set of functions with known
curvature (affine, convex, or concave) and monotonicity, and a composition rule
for combining the functions to produce expressions that are also convex or
concave. Some software does exist for solving quasiconvex problems (e.g.,
YALMIP \cite{yalmip}), but no DSLs exist for specifying them in a way that
guarantees quasiconvexity.

In this paper, we introduce \emph{disciplined quasiconvex programming}
(DQCP), an analog of DCP for quasiconvex optimization. Like DCP, DQCP is a
grammar that consists of a set of functions and rules for combining them. A
contribution of this paper is the development of a theorem for the composition
of a quasiconvex function with convex (and concave) functions that guarantees
quasiconvexity of the composition. This rule includes as a special case the
composition rule for convex functions upon which DCP is based. The class of
programs producible by DQCP is a subset of QCPs (and depends on the function
library), and a superset of the class corresponding to DCP\@.

In
\secref{s-qcvx}, we review properties of quasiconvex functions, state our
composition theorem, and provide several examples of quasiconvex functions. In
\secref{s-sol}, we describe a bisection method for solving QCPs\@. In
\secref{s-dqcp}, we present DQCP\@, and in \secref{s-impl}, we describe
an implementation of DQCP in CVXPY 1.0.

\section{Quasiconvexity}\label{s-qcvx}
\subsection{Properties}\label{s-prop}
In this section, we review basic properties of quasiconvex functions, many
of which are parallels of properties of convex functions; see
\cite{greenberg1971} for many more.
\paragraph{Jensen's inequality.} Quasiconvex functions are characterized by a
kind of Jensen's inequality: a function $f$ mapping a set $C$ into $\RR$ is
quasiconvex if and only if $C$ is convex and, for any $x, y \in C$ and $\theta
\in [0, 1]$,
\[
f(\theta x + (1 - \theta)y) \leq \max \{f(x), f(y)\}.
\]
Similarly, $f$ is quasiconcave if and only if $C$ is convex and
$f(\theta x + (1 - \theta)y) \geq \min \{f(x), f(y)\}$, for all $x, y \in C$
and $\theta \in [0, 1]$.

\paragraph{Functions on the real line.} For $f : C \subseteq \RR \to \RR$,
quasiconvexity can be described in simple terms: $f$ is quasiconvex if
it is nondecreasing, nonincreasing, or nonincreasing over $C \cap (\infty,
t]$ and nondecreasing over $[t, \infty) \cap C$, for some $t \in C$.

\paragraph{Representation via a family of convex functions.}
The sublevel sets of a quasiconvex function can be represented as inequalities
of convex functions. In this sense, every quasiconvex function can be represented
by a family of convex functions. If $f: C \to \RR$ is quasiconvex,
then there exists a family of convex functions $\phi_t : C \to \RR$, indexed by
$t \in \RR$, such that
\[
f(x) \leq t \iff \phi_t(x) \leq 0.
\]
The indicator functions for the sublevel sets of $f$,
\[
\phi_t(x) = \begin{cases}
0 & f(x) \leq t \\
\infty & \text{otherwise},
\end{cases}
\]
generate one such family. As another example, if the sublevel sets of $f$ are
closed, a suitable family is $\phi_t(x) = \inf_{z \in \Set{z}{f(z) \leq t}}
\norm{x - z}$. We are typically interested in finding families that possess
nice properties. For the purpose of DQCP, we seek functions $\phi_t$ whose
$0$-sublevel sets can be represented by convex cones over which optimization is
tractable.

\paragraph{Partial minimization.} Minimizing a quasiconvex function over
a convex set with respect to some of its variables yields another
quasiconvex function.

\paragraph{Supremum of quasiconvex functions.} The supremum of a family
of quasiconvex functions is quasiconvex, as can be easily verified
\cite[\S3.4.4]{boyd2004}; similarly, the infimum of quasiconcave functions is
quasiconcave.

\paragraph{Composition with monotone functions.} If $g : C \to \RR$ is
quasiconvex and $h$ is a nondecreasing real-valued function on the real line,
then $f = h \circ g$ is quasiconvex. This can be seen by observing that for any
$\alpha \in \RR$, a point $x$ (belonging to the domain of $f$) is in the
$\alpha$-sublevel set of $f$ if and only if
\[
g(x) \leq \sup \Set{y}{h(y) \leq \alpha}.
\]
Because $g$ is quasiconvex, this shows that the sublevel sets of $f$ are
convex. Similarly, a nonincreasing function of a quasiconvex function is is
quasiconcave, a nondecreasing function of a quasiconcave function is
quasiconcave, and a nonincreasing function of a quasiconcave function is
quasiconvex.

\subsection{Composition theorem}\label{s-comp}
A basic result from convex analysis is that
a nondecreasing convex function of a convex function is convex; DCP is based on
a generalization of this result. The composition rule for convex functions
admits a partial extension for quasiconvex functions, which we state below as a theorem.
Though the theorem is straightforward, we are unaware of any references to it
in the literature. Of course, the analog of the theorem for quasiconcave
functions also holds.

In the statement of the
theorem, when considering a function $g$ mapping a subset of $\RR^n$ into
$\RR^k$, we use $g_1, g_2, \ldots, g_k$ to denote the components of $g$. These
components are the real functions defined by
\[
g(x) = (g_1(x), g_2(x), \ldots, g_k(x)),
\]
for $x$ in the domain of $g$.

\begin{theorem}\label{thm-comp}
Suppose $h$ is a quasiconvex mapping of a subset $C$ of $\RR^k$ into $\RR \cup
\infty$, and $\{I_1, I_2, I_3\}$ is a partition of $\{1, 2, \ldots, k\}$ such
that $h$ is nondecreasing in the arguments indexed by $I_1$ and
nonincreasing in the arguments indexed by $I_2$. Suppose also that $g$ maps a subset
of $\RR^n$ into $\RR^k$ in such a way that its components $g_i$ are convex for $i \in
I_1$, concave for $i \in I_2$, and affine for $i \in I_3$. Then the composition
\[
f = h \circ g
\]
is quasiconvex. If additionally $h$ is convex, then $f$ is convex as well.
\end{theorem}
The final statement of the theorem is just the well-known composition rule for
convex functions.

We provide two proofs of this result. The first proof directly verifies that
the domain of $f$ is convex and that $f$ satisfies the modified Jensen's
inequality. This proof is almost identical to a proof of the composition
theorem for convex functions. The only difference is that an application of
Jensen's inequality for convex functions is replaced with its variant for
quasiconvex functions. The second proof just applies the composition theorem
for convex functions to the representation of a quasiconvex function via a
family of convex functions.

\begin{proof}[Proof via Jensen's inequality.]
Assume $x, y$ are in the domain of $f$, and
$\theta \in [0, 1]$. Since the components of $g$ are convex or concave (or
affine), the convex combination $\theta x + (1 - \theta)y$ is in the domain of $g$. For
$i \in I_1$, the components are convex, so
\[
g_i(\theta x + (1 - \theta)y) \leq \theta g_i(x) + (1 - \theta) g_i(y).
\]
For $i \in I_2$, the inequality is reversed, and for $i \in I_3$, it is an
equality. Since $x$ and $y$ are in the domain of $f$, $g(x)$ and $g(y)$ are in the
domain $C$ of $h$, and $\theta g(x) + (1
- \theta) g(y) \in C$. Let $e_i$ denote the $i$th standard basis vector of
  $\RR^k$. Since $h$ is an extended-value function and in light of its
per-argument monotonicities, $C$ extends infinitely in the directions $-e_i$ for
$i \in I_1$ and $e_i$ for $i \in I_2$. This fact, combined with the
inequalities involving the components of $g$ and the fact that $\theta g(x) +
(1 - \theta)g(y) \in C$, shows that $g(\theta x + (1 - \theta)y) \in C$. Hence
the domain of $f$ is convex.

By the monotonicity of $h$ and Jensen's inequality applied to the components of
$g$,
\[
h(g(\theta x + (1 - \theta)y)) \leq h(\theta g(x) + (1 - \theta) g(y)).
\]
Because $h$ is quasiconvex,
\[
h(\theta g(x) + (1 - \theta) g(y)) \leq \max \{h(g(x)), h(g(y))\}.
\]
Hence $f$ is quasiconvex.
\end{proof}

\begin{proof}[Proof via representation by convex functions.]
Let $\phi_t : C \to \RR$ be a member of a family of convex functions, indexed
by $t$, such that $\phi_t(x) \leq 0$ if and only if $h(x) \leq t$. Assume
without loss of generality that the per-argument monotonicities of $\phi_t$
match those of $h$ (\eg, take $\phi_t$ to be the indicator function for the
$t$-sublevel set of $h$). Then $f(x) = h(g(x)) \leq t$ if and only if
$\phi_t(g(x)) \leq 0$. By the composition theorem for convex functions,
$\phi_t \circ g$ is convex. We therefore conclude that the sublevel sets of $f$
are convex, \ie, $f$ is quasiconvex.
\end{proof}

\subsection{Examples}\label{s-ex}

\paragraph{Product.} The scalar product $f(x, y) = xy$ is quasiconcave when
restricted to either $\RR^2_{+}$ or $\RR^2_{-}$, where $\RR^n_{+}$ denotes the
set of nonnegative real $n$-vectors and $\RR^n_{-}$ the set of nonpositive
real $n$-vectors. The product is quasiconvex when one variable is nonnegative and
the other is nonpositive. From this fact and the composition rule, one can
deduce that the product of two nonnegative concave functions is quasiconcave
(see also \cite{bector1968, kantrowitz2005}), and the product of a
nonnegative concave function with a nonpositive convex function is quasiconvex.

\paragraph{Ratio.} The ratio $f(x, y) = x / y$ is quasilinear on $\RR \times \RR_{++}$,
as well as on $\RR \times \RR_{--}$ (but not on $\RR^2$), where $\RR^n_{++}$ and
$\RR^n_{--}$ denote the sets of positive and negative real $n$-vectors,
respectively. When $x \geq 0$ and $y > 0$, $f$ is increasing in $x$ and
decreasing in $y$. Hence the ratio of a nonnegative convex function and a
positive concave function is quasiconvex, and the ratio of a nonnegative
concave function and a positive convex function is quasiconcave. The problem of
maximizing the ratio of a nonnegative concave function and a positive convex
function is known as concave-fractional programming \cite{schaible1978,
schaible1981}.

\paragraph{Linear-fractional function.} The function
\[
f(x) = \frac{a^T x + b}{c^Tx + d}
\]
is quasilinear when the denominator is positive. This can be seen by the
composition rule, since the ratio $x/y$ is quasilinear when $y > 0$. It is also
quasilinear when restricted to negative denominators. The problem of minimizing
a linear-fractional function over a polyhedron is known as linear-fractional
programming. Though linear-fractional programming is often described as a
generalization of linear programming, linear-fractional programs can be reduced
to linear programs \cite{charnes1962}.

\paragraph{Distance ratio function.} The function
\[
f(x) = \frac{\norm{x - a}_2}{\norm{x - b}_2}
\]
is quasiconvex on the halfspace $\Set{x \in \RR^n}{\norm{x - a}_2 \leq \norm{x -
b}_2}$. This result cannot be derived by applying the composition rule to the
ratio function, but it is simple to show that its sublevel sets are
Euclidean balls \cite[\S3.4]{boyd2004}.

\paragraph{Monotone functions on the real line.} Monotone functions whose domains
are convex subsets of $\RR$ are quasilinear; examples include the exponential
function, logarithm, square root, and positive odd powers.

\paragraph{Generalized eigenvalue.} The maximum eigenvalue
of a symmetric matrix is convex, since it can be written as the supremum of
a family of linear functions. Analogously, the maximum generalized eigenvalue
$\lambda_{\mathrm{max}}(A, B)$ of a pair of symmetric matrices $(A, B)$ (with $B$
positive definite) is quasiconvex, since
\[
\lambda_{\mathrm{max}}(A, B)  = \sup_{x \neq 0} \frac{x^TAx}{x^TBx}
\]
is the supremum of a family of linear-fractional functions
\cite[\S3.4]{boyd2004}. Another way to see this is to note that the inequality
\[
\lambda_{\mathrm{max}}(A, B)  = \sup\{\lambda \in \RR \mid Ax = \lambda Bx\} \leq t
\]
is satisfied if and only if $tB - A$ is positive semidefinite. Similarly, the
minimum generalized eigenvector is quasiconcave in $A$ and $B$.

\subsubsection{Integer-valued functions}
\paragraph{Ceiling and floor.} The functions
$\ceil{x} = \inf \Set{z \in \ZZ}{z \geq x}$ and
$\floor{x} = \sup \Set{z \in \ZZ}{z \leq x}$ are quasilinear, because they are
monotone functions on the real line.

\paragraph{Sign.} The function mapping a real number to $-1$ if it is negative
and $+1$ otherwise is quasilinear.

\paragraph{Rectangle.} The rectangle function $f : \RR \to \RR$ given by
\[
f(x) = \begin{cases}
0 & |x| > \frac{1}{2} \\
1 & |x| \leq \frac{1}{2}
\end{cases}
\]
is quasiconcave.

\paragraph{Length of a vector.} The length of a vector in $\RR^n$ is defined as
the largest index corresponding to a nonzero component:
\[
\operatorname{len}(x) = \max \{i \mid x_i \neq 0\}.
\]
This function is quasiconvex on $\RR^n$ because its sublevel sets are
subspaces. The inequality $f(x) \leq \alpha$ implies
$x_i = 0$ for $i = \floor{\alpha} + 1, \ldots, n$.

\paragraph{Cardinality of a nonnegative vector.} The function
$\operatorname{card}(x)$, which gives the number of nonzero components in
the vector $x$, is quasiconcave on $\RR^n_+$: $\operatorname{card}(x + y) \geq
\min\{\operatorname{card}(x), \operatorname{card}(y)\}$ for nonnegative
$x$ and $y$.

\paragraph{Matrix rank.} The matrix rank is
quasiconcave on the set of positive semidefinite matrices, since the rank
of a sum of positive semidefinite matrices is at least the minimum of the ranks
of the matrices.

\section{Solution method}\label{s-sol}
The problem of minimizing a quasiconvex function $f : C \to \RR$ can be solved in many ways
\cite{kiwiel2001, konnov2003, hazan2015}. Here, we describe a simple method
that reduces a QCP to a
sequence of convex feasibility problems \cite[\S4.2.5]{boyd2004}. Suppose the
interval $[\alpha, \beta]$ is known to contain the optimal value $p^\star$.
Put $t = (\alpha + \beta)/2$, and let $\phi_t : C \to \RR$ be a family of convex
functions indexed by $t \in \RR$ such that $f(x) \leq t$ if and only if
$\phi_t(x) \leq 0$. Consider the convex feasibility problem
\begin{equation}\label{eqn-feas}
\begin{array}{ll}
\mbox{find} & x \\
\mbox{subject to} & \phi_t(x) \leq 0. \\
\end{array}
\end{equation}
If this problem yields a feasible point $x$, then $p^\star \leq t$ and in
particular $p^\star \in [\alpha, f(x)]$; otherwise, $p^\star \in [t, \beta]$.
In either case, solving the feasibility problem yields an interval containing
the optimal value, with width half as large as the original interval. To obtain
an $\epsilon$-suboptimal solution to the QCP, we repeat this process until the
width of the interval is at most $\epsilon$, which requires at most
$\ceil{\log_2(\beta - \alpha)/\epsilon}$ iterations.

\paragraph{Finding an initial interval for bisection.}
The optimal value $p^\star$ is usually not known before solving a QCP\@. In such
cases, a simple heuristic can be employed to find an interval containing it,
assuming that the QCP is feasible (which can be checked by solving a single
convex feasibility
problem). Start with a candidate interval $[\alpha, \beta]$, where $\alpha <
0$ and $\beta > 0$. If the
problem $(\ref{eqn-feas})$ is feasible for $t=\beta$ and infeasible for
$t=\alpha$, then $p^\star \in [\alpha, \beta]$. Otherwise, if the
problem is infeasible for $t=\beta$, put $\alpha := \beta$ and $\beta := 2\beta
$. If on the other hand the problem is feasible for $t=\alpha$, put $\beta
:= \alpha$ and $\alpha := 2\alpha$. Repeating this process will eventually
produce an interval containing $p^\star$, provided that the QCP is not
unbounded.

\section{Disciplined quasiconvex programming}\label{s-dqcp}

DQCP is a grammar for constructing QCPs
from a set of functions, or \emph{atoms}, with known curvature (affine, convex,
concave, quasiconvex, or quasiconcave) and per-argument monotonicities. A
program produced using DQCP is called a disciplined quasiconvex program; we say
that such programs are DQCP-compliant, or just DQCP, for short. DQCP guarantees
that every function appearing in a disciplined quasiconvex program is affine,
convex, concave, quasiconvex, or quasiconcave.

A disciplined quasiconvex program is an optimization problem of the form
\begin{equation}\label{eqn-dqcp}
\begin{array}{ll}
\mbox{minimize} & f_0(x) \\
\mbox{subject to} & f_i(x) \leq \alpha_i, \quad i=1, \ldots, m_1\\
 & \beta_i \leq g_i(x), \quad i=1, \ldots, m_2\\
 & \tilde{f}_i(x) \leq \tilde{g}_i(x), \quad i=1, \ldots, m_3\\
& h_i(x) =\tilde{h}_i(x), \quad i=1, \ldots, p.
\end{array}
\end{equation}
The functions $f_i$ must be quasiconvex, $g_i$ must be quasiconcave, $\tilde{f}_i$
must be convex, $\tilde{g}_i$ must be concave, and $h_i$, $\tilde{h}_i$
must be affine; $\alpha_i$ and $\beta_i$ must be constants. All of the
functions appearing in (\ref{eqn-dqcp}) must be produced from atoms, using
only the composition rule from Theorem~\ref{thm-comp} and the rules governing
the maximum of quasiconvex functions, the minimum of quasiconcave functions,
and composition with monotone functions (see \S\ref{s-prop}). Because
Theorem~\ref{thm-comp} includes the composition rule for convex functions as a
special case, DQCP is a modest extension of DCP.

A mathematical expression is verifiably quasiconvex under DQCP if it is
\begin{itemize}
\item a convex expression;
\item a quasiconvex atom, applied to a variable or constant;
\item the max of quasiconvex expressions;
\item a nondecreasing function of a quasiconvex expression, or a nonincreasing
      function of a quasiconcave expression;
\item the composition of a quasiconvex atom with convex, concave, and
      affine expressions that satisfies the hypotheses of Theorem~\ref{thm-comp}.
\end{itemize}
These rules are applied recursively, with the recursion bottoming out
at variables and constants. For example, if $\exp(\cdot)$
and the generalized eigenvalue $\lmbdmax(\cdot, \cdot)$ are atoms, and $X$ and $Y$
are matrix variables, then the expressions
\[
\lmbdmax(X, Y), \quad \exp(\lmbdmax(X, Y)), \quad \textnormal{and} \quad \exp(\exp(\lmbdmax(X, Y)))
\]
are all verifiably quasiconvex under DQCP, since $\exp(\cdot)$ is increasing
and $\lmbdmax(\cdot, \cdot)$ is quasiconvex. Likewise, an expression is
quasiconcave under DQCP if it is a concave expression, a quasiconcave atom
applied to a variable or constant, the min of quasiconcave functions, a
nondecreasing function of a quasiconcave function, a nonincreasing function of
a quasiconvex function, or a valid composition of a quasiconcave function
with convex, concave, and affine functions. Whether an expression is convex,
concave, or affine under DQCP is precisely the same as under DCP.

\begin{figure}
\adjustbox{valign=t}{\begin{minipage}[c]{0.40\textwidth}
\includegraphics[height=0.38\textheight]{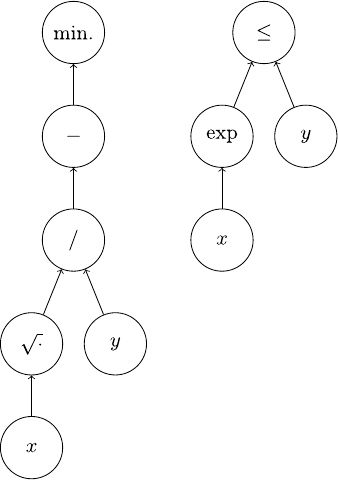}
\end{minipage}}%
\hfill
\adjustbox{valign=t}{\begin{minipage}[b]{0.60\textwidth}
\caption{Expression trees representing the program
(\ref{eqn-simple-ex}).
}
\label{fig-dqcp}
\end{minipage}}
\end{figure}

A DQCP program is naturally represented as a collection of expression trees,
one for the objective and one for each constraint. Verifying whether a
program is DCQP amounts to recursively verifying that each expression tree
is DQCP\@. For example, the program
\begin{equation}\label{eqn-simple-ex}
\begin{array}{ll}
\mbox{minimize} & -\sqrt{x}/y \\
\mbox{subject to} &  \exp(x) \leq y
\end{array}
\end{equation}
can be represented by the trees shown in figure~\ref{fig-dqcp}. This program is
DQCP when $y$ is known to be positive, because the ratio of a nonnegative
concave functions and a positive convex function is quasiconcave, and the
negation of a quasiconcave function is quasiconvex. The atoms in this program
are the functions $\exp(\cdot)$, $\sqrt{\cdot}$, and $\cdot/\cdot$.

Every disciplined quasiconvex program is a QCP, but the converse is not
true. This is not a limitation in practice, since the atom library is
extensible.

\paragraph{The grammar.} Table~\ref{tbl-grm} specifies the DQCP grammar, in the
programming languages sense \cite[\S4]{dragonbook}. In
the specification, $\mathrm{S}$ denotes the start symbol. The symbols
\[
\AFF,\quad \CVX,\quad \CCV,\quad \QCVX,\quad \QCCV
\]
are nonterminals used to represent affine, convex, concave, quasiconvex, and
quasiconcave expressions producible by DQCP. Their lowercase
counterparts represent atoms, \eg, $\mathrm{cvx}$ stands for
a convex atom.  Atoms can have multiple curvatures.
For example, every affine atom is also a convex atom and a concave atom. The
symbols
\[
\mathrm{incr}, \quad \mathrm{decr}
\]
denote nondecreasing and nonincreasing functions, respectively,
\[
\mathrm{constant}, \quad \mathrm{variable}
\]
denote numerical constants and optimization variables, and
\[
\mathrm{cvx}(\CVX, \ldots, \CVX, \CCV, \ldots, \CCV, \AFF, \ldots, \AFF)
\]
denotes a composition of a convex atom with convex,
concave, and affine expressions that can be certified as convex via
Theorem~\ref{thm-comp}. Because DQCP is a grammar for QCPs, it can
be used to define the syntax of a DSL for quasiconvex optimization.
\begin{table}[]
\begin{tabular}{ll}
$\mathrm{S}$ & $\rightarrow$ $\QCVX$ \\
$\mathrm{S}$ & $\rightarrow$ $\QCCV$ \\ \\
$\LEAF$ & $\rightarrow$ $\mathrm{constant}$ \\
$\LEAF$ & $\rightarrow$ $\mathrm{variable}$ \\ \\
$\AFF$ & $\rightarrow$ $\LEAF$ \\
$\AFF$ & $\rightarrow$ $\mathrm{aff}(\AFF, \ldots, \AFF)$ \\ \\
$\CVX$ & $\rightarrow \AFF$ \\
$\CVX$ & $\rightarrow  \mathrm{cvx}(\CVX, \ldots, \CVX, \CCV, \ldots, \CCV, \AFF, \ldots, \AFF)$ \\ \\
$\CCV$ & $\rightarrow \AFF$ \\
$\CCV$ & $\rightarrow  \mathrm{ccv}(\CCV, \ldots, \CCV, \CVX, \ldots, \CVX, \AFF, \ldots, \AFF)$ \\ \\
$\QCVX$ & $\rightarrow \CVX$ \\
$\QCVX$ & $\rightarrow  \mathrm{qcvx}(\CVX, \ldots, \CVX, \CCV, \ldots, \CCV, \AFF, \ldots, \AFF)$\\
$\QCVX$ & $\rightarrow  \mathrm{incr}(\QCVX)$ \\
$\QCVX$ & $\rightarrow  \mathrm{decr}(\QCCV)$ \\
$\QCVX$ & $\rightarrow  \max\{\QCVX, \ldots, \QCVX\}$ \\ \\
$\QCCV$ & $\rightarrow \CCV$ \\
$\QCCV$ & $\rightarrow \mathrm{qccv}(\CCV, \ldots, \CCV, \CVX, \ldots, \CVX, \AFF, \ldots, \AFF)$ \\
$\QCCV$ & $\rightarrow  \mathrm{incr}(\QCCV)$ \\
$\QCCV$ & $\rightarrow  \mathrm{decr}(\QCVX)$ \\
$\QCCV$ & $\rightarrow  \min\{\QCCV, \ldots, \QCCV\}$
\end{tabular}
\caption{The DQCP grammar, which extends DCP\@. The rules for
compositions with convex, concave, and affine expressions denote compositions
satisfying the hypotheses of Theorem~\ref{thm-comp}.} \label{tbl-grm}
\end{table}

\section{Implementation}\label{s-impl}
We have implemented DQCP in CVXPY 1.0, a Python-embedded DSL for convex
optimization \cite{cvxpy, cvxpyrewriting}. Our implementation, which is
available at
\begin{center}
  \url{https://www.cvxpy.org},
\end{center}
makes CVXPY the first DSL for quasiconvex optimization. Because DQCP is a
generalization of DCP, it fits seamlessly into CVXPY, which parses
problems using DCP by default. Our atom library includes many of the functions
presented in \secref{s-ex}. We have also implemented the bisection method
described in \secref{s-sol}.

\subsection{Canonicalization} The process of rewriting a problem to an
equivalent standard form is called canonicalization. In CVXPY 1.0,
canonicalization is facilitated by \texttt{Reduction} objects, which
rewrite problems of one form into equivalent problems of another form.

We have implemented a reduction called \texttt{Dqcp2Dcp} that canonicalizes
DQCP problems by converting them into an equivalent one-parameter
family of DCP feasibility problems.
When applied to a DQCP problem, this reduction first introduces a scalar
parameter and constrains the problem's objective to be no greater than the
parameter. It recursively processes this constraint and every other
constraint, representing the sublevel sets of quasiconvex expressions and
superlevel sets of quasiconcave expressions in DCP-complaint ways. The reduction
then emits a parameterized DCP problem. The constraints of the emitted problem are
the canonicalized constraints of the original problem, and the objective is to
find an assignment to the variables that satisfies the constraints. A solution
to the original problem can be obtained by running bisection on the emitted
problem.

\subsection{Bisection} We have implemented the bisection routine described
in \S\ref{s-sol}. Our method first checks whether the original problem is
feasible by solving a convex feasibility problem. If the problem is feasible,
our routine automatically finds an interval containing the optimal value and then runs
bisection. Our bisection routine tightens the boundaries of the bisection
interval depending on the values of the original problem's objective function.
For example, when the objective is integer-valued, our implementation will
tighten a lower bound $\alpha$ to $\ceil{\alpha}$, and an upper bound $\beta$
to $\floor{\beta}$.

\subsection{Examples}
\paragraph{Hello, world.} Below is an example of how to use CVXPY 1.0 to specify
and solve the problem~(\ref{eqn-simple-ex}), meant to highlight the syntax of our
modeling language. More interesting examples are subsequently presented.

\begin{lstlisting}[xleftmargin=.05\textwidth, xrightmargin=.2\textwidth]
import cvxpy as cp

x = cp.Variable()
y = cp.Variable(pos=True)
objective_fn = -cp.sqrt(x)/y
objective = cp.Minimize(objective_fn)
constraint = cp.exp(x) <= y
problem = cp.Problem(objective, [constraint])
problem.solve(qcp=True)
print("Optimal value: ", problem.value)
print("x: ", x.value)
print("y: ", y.value)
\end{lstlisting}
The optimization problem \texttt{problem} has two scalar variables, x and y.
Notice that \texttt{y} is declared as positive in line \texttt{3}, with
\texttt{pos=True}. The objective is to minimize the ratio of $-\sqrt{x}$ and $y$,
which is quasiconvex since the ratio is
quasiconcave when the numerator is a nonnegative concave expression and the
denominator is a positive convex expression. Line \texttt{6} constructs
the objective of the problem. In
line \texttt{7}, \texttt{exp(x)} is constrained to be no larger than
\texttt{y} via the relational operator \texttt{<=}. Line 8 constructs
\texttt{problem}, which represents the optimization problem as two expression
trees, one for \texttt{objective\_fn} and one for \texttt{constraint}. The
internal nodes in these expression trees are the atoms \texttt{sqrt},
\texttt{exp}, ratio (\texttt{/}), and negation $(\texttt{-})$. The problem
is DQCP, which can be verified by asserting \texttt{problem.is\_dqcp()}.
Line \texttt{9} canonicalizes \texttt{problem}, parsing
it as a DQCP (\texttt{qcp=True}), and then solves it by bisection. The optimal
value of the problem and the values of the variables are printed in lines
10-12, yielding the following output.
\begin{lstlisting}[xleftmargin=.05\textwidth, xrightmargin=.2\textwidth]
Optimal value:  -0.4288821220397949
x:  0.49999737143004713
y:  1.648717724845007
\end{lstlisting}

As this example makes clear, users do not need to know how canonicalization
or bisection work. All they need to know is how to construct DQCP problems.
Calling the \texttt{solve} method on a \texttt{Problem} instance with the
keyword argument \texttt{qcp=True} canonicalizes the problem and retrieves a
solution. If the user forgets to type \texttt{qcp=True} when her problem is
DQCP (and not DCP), a helpful error message is raised to alert her of the
omission.

\paragraph{Generalized eigenvalue matrix completion.} We have implemented the
maximum generalized eigenvalue as an atom. As an example, we can use CVXPY 1.0
to formulate and solve a generalized eigenvalue matrix completion problem. In
this problem, we are given \textit{some} entries of two symmetric matrices
$A$ and $B$, and the goal is to choose the missing entries so as to minimize
the maximum generalized eigenvalue $\lambda_{\mathrm{max}}(A, B)$. Letting $\Omega$
denote the set of indices $(i, j)$ for which $A_{ij}$ and $B_{ij}$ are known,
the optimization problem is
\begin{equation*}
\begin{array}{ll}
\mbox{minimize} & \lambda_{\mathrm{max}}(X, Y) \\
\mbox{subject to} & X_{ij} = A_{ij}, \, (i, j) \in \Omega, \\
                  & Y_{ij} = B_{ij}, \, (i, j) \in \Omega,
\end{array}
\end{equation*}
which is a QCP.  Below is an implementation of
this problem, with specific problem data
\begin{equation*}
\centering
A = \begin{bmatrix}
1.0 & ? &  1.9 \\
? & 0.8 &  ? \\
? & ?  &  ?
\end{bmatrix}, \quad
B = \begin{bmatrix}
3.4 & ? &  1.4 \\
? & 0.2 &  ? \\
? & ?  &  ?
\end{bmatrix}.
\end{equation*}
(The question marks denote the missing entries.)

\begin{lstlisting}[xleftmargin=.05\textwidth, xrightmargin=.2\textwidth]
import cvxpy as cp

X = cp.Variable((3, 3))
Y = cp.Variable((3, 3))
gen_lambda_max = cp.gen_lambda_max(X, Y)
omega = tuple(zip(*[[0, 0], [0, 2], [1, 1]]))
constraints = [
  X[omega] == [1.0, 1.9, 0.8],
  Y[omega] == [3.0, 1.4, 0.2],
]
problem = cp.Problem(cp.Minimize(gen_lambda_max), constraints)
problem.solve(qcp=True)
print("Generalized eigenvalue: ", gen_lambda_max.value)
print("X: ", X.value)
print("Y: ", Y.value)
\end{lstlisting}
Executing the above code prints the below output.
\begin{lstlisting}[xleftmargin=.05\textwidth, xrightmargin=.2\textwidth,
numbers=none]
Objective:  4.000002716411653
X:  [[9.99999767e-01 9.86154616e-16 1.89999959e+00]
 [9.86154616e-16 7.99999761e-01 5.19126535e-15]
 [1.89999911e+00 5.19126535e-15 1.25733692e+00]]
Y:  [[ 2.99999980e+00 -2.78810135e-16  1.40000015e+00]
 [-2.78810135e-16  1.99999804e-01  2.14473098e-16]
 [ 1.40000015e+00  2.14473098e-16  1.14038551e+00]]
\end{lstlisting}
Notice that the \texttt{gen\_lambda\_max} atom automatically enforced the
symmetry and positive definiteness constraints on $X$ and $Y$.

\paragraph{Minimum length least squares.} Our atom library includes several
integer-valued functions, including the length function. As an example, the
following QCP finds a minimum-length vector $x \in \RR^n$ that has small
mean-square error for a particular least squares problem:
\begin{equation*}
\begin{array}{ll}
\mbox{minimize} & \operatorname{len}(x) \\
\mbox{subject to} & 1/n\norm{Ax - b}_2^2 \leq \epsilon.
\end{array}
\end{equation*}
The problem data are $A \in \RR^{n \times n}$, $b \in \RR^n$, and $\epsilon \in
\RR$. Below is an implementation of this problem in CVXPY.
\begin{lstlisting}[xleftmargin=.05\textwidth, xrightmargin=.2\textwidth]
import cvxpy as cp
import numpy as np
np.set_printoptions(precision=2)

n = 10
np.random.seed(1)
A = np.random.randn(n, n)
x_star = np.random.randn(n)
b = A @ x_star
epsilon = 1e-2

x = cp.Variable(n)
mse = cp.sum_squares(A @ x - b)/n
problem = cp.Problem(cp.Minimize(cp.length(x)), [mse <= epsilon])
problem.solve(qcp=True)
print("Length of x: ", problem.value)
print("MSE: ", mse.value)
print("x: ", x.value)
print("x_star: ", x_star)
\end{lstlisting}

Running the code produces the following output.
\begin{lstlisting}[xleftmargin=.05\textwidth, xrightmargin=.2\textwidth, numbers=none]
Length of x:  8.0
MSE:  0.00926009328775564
x:  [-0.26  1.38  0.21  0.94 -1.15  0.15  0.66 -1.16 -0.   -0.  ]
x_star:  [-0.45  1.22  0.4   0.59 -1.09  0.17  0.74 -0.95 -0.27  0.03]
\end{lstlisting}

\section*{Acknowledgments} The authors thank Steven Diamond for
many useful discussions.

\printbibliography{}

\end{document}